\documentclass[12pt,reqno]{amsproc}

\newtheorem{theo}{Theorem}
\newtheorem{rem}{Remark}

\newtheorem{lem}{Lemma}

\newtheorem{df}{Definition}

\newcommand\eps\varepsilon
\newcommand\ph\varphi
\newcommand\kap\Lambda

\usepackage{enumerate}
\usepackage{graphicx}
\usepackage{float}

\begin{document}

\title[Monotone Discontinuous ODE in Sequence Spaces]
{Monotone  ODEs with Discontinuous Vector Fields in Sequence Spaces }

\author[Oleg Zubelevich]{Oleg Zubelevich\\
 Steklov Mathematical Institute of Russian Academy of Sciences \\oezubel@gmail.com
 }
\email{oezubel@gmail.com}
\date{}
\thanks{The research was funded by a grant from the Russian Science
Foundation (Project No. 19-71-30012)}

\begin{abstract}We consider a system of ODE in a Fr\'echet space with unconditional Schauder basis. The right side of the ODE is a discontinuous function. Under certain monotonicity conditions we prove an existence theorem for the corresponding initial value problem. We employ an idea of the partial order which seems to be new in this field.
\end{abstract}

\maketitle
\numberwithin{equation}{section}
\newtheorem{theorem}{Theorem}[section]
\newtheorem{lemma}[theorem]{Lemma}
\newtheorem{definition}{Definition}[section]

\section{Introduction}The analysis of ODE with non Lipschitz right hand side has a long history.
Without any claims on  completeness of exposition we just note some principle points of this history. A detailed discussion  of  further developments in any  of these points requires a separate survey.

The first result belongs to G. Peano (1890).  G. Peano considered an initial value problem
\begin{equation}\label{sxfg00gg}\dot x=f(t,x),\quad x(t_0)=x_0\end{equation} where $f$ is a continuous mapping of some domain $$D\subset\mathbb{R}^{m+1}=\{(t,x)\},\quad x=(x_1,\ldots,x_m)$$ with values in $\mathbb{R}^{m}$.

G. Peano stated that this problem has a solution  that is defined locally for small $|t-t_0|$. This solution may not be unique.

C. Carath\'eodory relaxed the conditions of this theorem up to  measurability of the function $f$ in $t$.

All these results are essentially  based on the fact that a closed ball in $\mathbb{R}^m$ is compact. In an infinite dimensional Banach space  they are in general invalid.

The corresponding example was first constructed  by J. Dieudonn\'e \cite{ded}.
His example is as follows:
$$\dot x_k=\sqrt{|x_k|}+\frac{1}{k},\quad x_k(0)=0,\quad k\in\mathbb{N},\quad t\ge 0.$$ By direct integration one  can easily  check that this IVP does not have solutions $x(t)=\{x_k(t)\}\in c_0$.

In the infinite dimensional case  besides the continuity    one must impose some extra
compactness conditions on $f$ \cite{mill}, \cite{szufla} in order to recover the existence.

Observe that  all the existence results mentioned above   follow in  one way or another from the Schauder-Tychonoff fixed point theorem.

Note also that in case of  non-Lipschitz equations  a solution  is in general not unique and one has  to make separate efforts to study the problem of uniqueness.

When   the function $f$ is discontinuous in $x$ then even examples that are  pretty innocuous from the first glance can
provide nonexistence.
Indeed, consider a scalar IVP
\begin{equation}\label{sfvvv}
\dot x=h(x),\quad x(0)=1,\end{equation}
where
$$h(x)=\begin{cases}
  1,  & \mbox{if } x\le 1, \\
  -1, & \mbox{if } x>1.
\end{cases}$$
This problem does not have  a continuous solution $x(t)$ in the sense of integral equation:
$$x(t)=1+\int_0^th(x(\xi))d\xi,\quad  t>0.$$
To show this assume the converse:  this solution exists for some $t>0$. It is clear that it can not be equal to 1 identically. Thus for some $t'>0$ we have $x(t')>1$. (The case $x(t')<1$ is accomplished  similarly.) Then  there exists an interval $(t_1,t')$ such that $$t\in (t_1,t')\Longrightarrow x(t)>1$$ and $x(t_1)=1.$

For $t\in (t_1,t')$ we can write
$$x(t)=x(t_1)+\int_{t_1}^tf(x(\xi))d\xi=1-(t-t_1)<1.$$ This is a contradiction.

Such examples prompt an idea to change   the concept of a solution itself.
Note in addition that if the right side of equation (\ref{sxfg00gg}) is just a measurable function then even for continuous $x(t)$ a mapping $t\mapsto f(t,x(t))$ is not obliged to be measurable \cite{folland}.

The corresponding transformation of the notion of a solution  was proposed by   A. Filippov \cite{filipp}. According to him an absolutely continuous function $x(t)$ is a solution to  (\ref{sxfg00gg}) with $f$ measurable if the following inclusion holds
$$\dot x(t)\in\bigcap_{r>0}\bigcap_{N}\mathrm{conv}\,f\big(t,B_r(x(t))\backslash N\big)\quad\mbox{ for almost all $t$. }$$
Here $B_r(x)\subset\mathbb{R}^m$ is an open ball of the radius $r$ and the center at $x$.   The intersection $\bigcap_{N}$ is taken over all measure-null sets $N$; and $\mathrm{conv}$ stands for the closed convex hull of a set.

Once we have denied  the classical concept of a solution then  a lot  of  reasonable generalizations arise.   Filippov's concept is good for control and for dry friction mechanics \cite{utlin}, \cite{zu_math_app}. A very different approach by DiPerna and Lions is good for PDE and fluid mechanics \cite{dip}.

In this article we try to save the classical concept of a solution  for some class of discontinuous ODE.

In the proofs we essentially use  the Lebesgue integration theory for functions that take values in Fr\'echet spaces.

The  Lebesgue integration theory of functions with values in Banach spaces is developed in \cite{lang}. The construction of \cite{lang} can easily be generalized to the case of Fr\'echet space valued functions.
Since we have not seen such a generalization   elsewhere and  for completeness  sake we collect a  necessary theory in section  \ref{xfhb77}.
There by means of minor modifications we adapt  the
theory from \cite{lang} to the case  of  Fr\'echet spaces.

\section{The  Main Theorems}

Let $E$ stand for a  Fr\'echet space. Its topology is defined by a collection of seminorms $\{\|\cdot\|_n\}_{ n\in\mathbb{N}}$. Recall that a Fr\'echet space is Hausdorff: for any $x\ne 0$ there exists $i$ such that $\|x\|_i\ne 0$. Moreover a topology of a Fr\'echet space is completely metrizable by the following metrics
$$\rho(x,y)=\sum_{k=1}^\infty \frac{1}{2^k}\min\{1,\|x-y\|_k\}.$$
Assume that the space $E$  possesses an unconditional Schauder basis $\{e_k\}_{k\in\mathbb{N}}$. Recall several definitions.
\begin{df}A sequence $\{e_k\}_{k\in\mathbb{N}}\subset E$ is called a Schauder basis  if for every $x\in E$ there is a unique sequence of  scalars $\{x_k\}_{k\in\mathbb{N}}$ such that
\begin{equation}\label{cdvfghyu8765fghjknhbgv}x=\sum_{k=1}^\infty x_ke_k.\end{equation}
This series is convergent in the topology of $E$.

We shall say that  $\{e_k\}_{k\in\mathbb{N}}$ is an unconditional basis if for any $x\in E$ and for any permutation $\pi:\mathbb{N}\to\mathbb{N}$ the sum
$$ \sum_{k=1}^\infty x_{\pi(k)}e_{\pi(k)}$$ is convergent to the same element.\end{df}
Introduce linear functions $e^j:E\to \mathbb{R}$ by the formula $e^j (x)=x_j.$ These functions are continuous \cite{edvards}.

We use $I_T$ to denote $[0,T],\quad T>0$.

Equip the space $E$ with a partial order $\ll$ as follows
$$x=\sum_{k=1}^\infty x_ke_k\ll y=\sum_{k=1}^\infty y_ke_k\Longleftrightarrow x_i\le y_i,\quad i\in\mathbb{N}.$$
\begin{df}We shall say that a function $g:E\to \mathbb{R}$ is left continuous if for all $ x'\in E$ and for all sequences $$x_k\to x',\quad x_k\ll x_{k+1},\quad k\in\mathbb{N}$$ one has
$$\lim_{k\to \infty}g(x_k)=g(x').$$\end{df}

The main object of our study is the following initial value problem
\begin{equation}\label{sdgf5500}
\dot x(t)=f\big(t,x(t)\big)=\sum_{k=1}^\infty f_k\big(t,x(t)\big)e_k,\quad x(0)=\hat x\in E.\end{equation}
Here the function  $f:I_T\times E\to E$ is  such that all the functions $f_k:I_T\times E\to\mathbb{R}$ are   left continuous in the second argument when the first one is fixed.

For any fixed $x$ the function $t\mapsto f(t,x)$ is integrable on $I_T$.

Assume that there exists an element $C=\sum_{k=1}^\infty C_ke_k\in E$ such that for any $(t,x)\in I_T\times E$
the following inequality  holds
\begin{equation}\label{fgb00}
f(t,x)\ll C.\end{equation}

Assume that there exists an element $x_*\in E$  that satisfies    the  inequality
\begin{equation}\label{dsfg22ww} x_*\ll \hat x+\int_0^tf(\xi,x_*)d\xi,\quad t\in I_T.\end{equation}

\begin{df}We shall say that a function $x\in C(I_T,E)$ is a solution to IVP (\ref{sdgf5500}) if a function
$t\mapsto f(t,x(t))$ is integrable in $I_T$ and  the following equation
$$
x(t)=\hat x+\int_0^tf\big(\xi,x(\xi)\big)d\xi,\quad t\in I_T$$ is satisfied.
\end{df}
\begin{theo}\label{dsfg5}In addition to the hypotheses above assume  that $f$ is monotone:
\begin{equation}\label{sxdfgiii}
x\ll y\Longrightarrow f(t,x)\ll f(t,y),\quad \forall x,y\in E, \quad \forall t\in I_T.\end{equation}
Then problem (\ref{sdgf5500}) has a solution  $x(t)$.\end{theo}
Theorem \ref{dsfg5} is proved in section \ref{sxdfg55}.

Condition (\ref{sxdfgiii}) is essential and can not be dropped: see example (\ref{sfvvv}).

To show that  condition (\ref{fgb00}) also matters we modify  Dieudonn\'e's example.

The IVP
$$\dot x_k=q(x_k)+\frac{1}{k},\quad x_k(0)=0,\quad k\in\mathbb{N},\quad t\ge 0$$ does not have solutions in $c_0$;
here the function
$$q(\xi)=\begin{cases}
  \sqrt \xi,  & \mbox{if } \xi\ge 0, \\
  0, & \mbox{if } \xi<0
\end{cases}$$
is monotone.

Let $V\subset E$ be a nonempty set. We formally put
$$\sup V:=\sum_{k=1}^\infty v_k^* e_k,\quad v_k^*=\sup\{e^k(v)\mid v\in V\}.$$
Let $Q\subset C(I_T,E)$ stand for a set of solutions to problem (\ref{sdgf5500}). From theorem \ref{dsfg5} we know that $Q\ne\emptyset$.

For each $t\in I_T$ introduce a set $Q(t)=\{v(t)\mid v\in Q\}\subset E.$
\begin{theo}\label{sdfrg55oo}Assume that the conditions of theorem \ref{dsfg5} are fulfilled.

Then for each $t\in I_T$ an element
$$\overline x(t)=\sup Q(t)\in E$$ is well defined and $\overline x$  is a solution to problem (\ref{sdgf5500}).\end{theo}

\begin{theo}\label{sxdfg1100}Assume in addition that $E$ is a reflexive space. Then for almost all $t$ the solution $x(t)$ from theorem \ref{dsfg5} is weakly  differentiable: there exists $\dot x(t)\in E$ such that
for   almost all $t$ one has
$$\Big(\psi,\frac{x(t+h)-x(t)}{h}\Big)\to (\psi,\dot x(t)),\quad h\to 0,\quad h\ne 0,\quad \forall\psi\in E'$$ and
$$\dot x(t)=f(t,x(t)).$$\end{theo}
Theorem \ref{sxdfg1100} is proved in section \ref{xdfb7khkhkh}.

\begin{rem}Theorems \ref{dsfg5}, \ref{sdfrg55oo}, \ref{sxdfg1100} remain valid if the functions $f_k$ are right continuous in the second argument.
\end{rem}

\subsection{An Example}
Introduce a function
$$H(\eta)=\begin{cases}
  -1,  & \mbox{if } \eta\le 0, \\
  1, & \mbox{if } \eta>0.
\end{cases}$$
Let  functions
$$n:\mathbb{Z}_+\to \mathbb{Z}_+,\quad \rho_k\in L^1(0,\infty),\quad k\in\mathbb{Z}_+=\mathbb{N}\cup\{0\}$$
be given.

Let $B=\{z\in\mathbb{C}\mid|z|<1\}$ stand for an open ball.  Let $E$ be a space of analytic functions
$$u:B\to\mathbb{C},\quad u(z)=\sum_{k=0}^\infty u_kz^k,\quad u_k\in\mathbb{R}$$ with a collection of seminorms
$$\|u\|_n=\max\{|u(z)|\,\mid\,|z|\le 1-1/n\},\quad n\in\mathbb{N}.$$
A set $\{z^k\},\quad k\in\mathbb{Z}_+$ is an unconditional   Schauder basis in $E$.

Consider the following initial value problem
\begin{equation}\label{dsfg66}
 u_t=\sum_{k=0}^\infty (k+1)^pH\big(u_{n(k)}+\rho_k(t)\big)z^k,\quad u\mid_{t=0}=\hat u(z),\end{equation}
where  the constant $p$ belongs to $\mathbb{N}$.

It is not hard to show that for any $\hat u\in E$ and $T>0$   problem (\ref{dsfg66}) satisfies  all the conditions of theorem \ref{dsfg5} and   has a solution
$$u(t,z)=\sum_{k=0}^\infty u_k(t)z^k$$
from $C(I_T,E)$. Particularly to satisfy condition (\ref{fgb00}) one should take
$$C=\sum_{k=0}^\infty(k+1)^pz^k\in E.$$

 \section{Proofs of the Theorems}We denote all inessential positive constants by the same letter $c$.
\subsection{Auxiliary Facts}The following theorem is essentially based on the assumption that the Schauder basis is unconditional. This theorem generalizes the corresponding result for Banach spaces \cite{lin}.
\begin{theo}[\cite{zu}]\label{e5b5j}Fix a sequence
$\lambda=\{\lambda_j\}_{j\in\mathbb{N}}\in \ell_\infty$.
Then  $$\mathcal M_\lambda x=\sum_{k=1}^\infty \lambda_k x_k e_k,\quad x=\sum_{k=1}^\infty x_k e_k\in E$$
is a bounded linear operator of $E$ to $E$ and for any number $i'$ there exists a number i and a positive constant $c$ both independent on $\lambda$ such that
$$\|\mathcal M_\lambda x\|_{i'}\le c\|\lambda\|_{\ell_\infty}\cdot\|x\|_i,\quad \forall x\in E.$$
\end{theo}
Particularly this theorem implies
$$y=\sum_{k=1}^\infty y_ke_k\in E\Longrightarrow |y|:=\sum_{k=1}^\infty|y_k|e_k\in E.$$
Other consequence of this theorem is as follows.
\begin{lem}\label{srg3ee}Assume that constant vectors $$a=\sum_{k=1}^\infty a_ke_k,\quad b=\sum_{k=1}^\infty b_ke_k\in E$$ are such that $a\ll b$. Then for any sequence of reals $y_k,\quad a_k\le y_k\le b_k$ an element
$$y=\sum_{k=1}^\infty y_ke_k\in E$$ is well defined. \end{lem}Indeed,
$$y=\mathcal M_\alpha a+\mathcal M_\beta b,$$ where
$$\alpha=\{\alpha_k\},\quad \beta=\{\beta_k\}\in \ell_\infty,\quad \alpha_k+\beta_k=1,\quad \alpha_k,\beta_k\ge 0.$$

\begin{lem}\label{dfg55}
Assume that constant vectors $$a=\sum_{k=1}^\infty a_ke_k,\quad b=\sum_{k=1}^\infty b_ke_k\in E$$ are such that $a\ll b$. Then an interval
$$[a,b]:=\{x\in E\mid a\ll x\ll b\}$$ is a compact set.

Moreover for any $i'\in\mathbb{N}$ there exists $i\in\mathbb{N}$ and a positive constant $c>0$ such that
\begin{equation}\label{dfgh77}x\in[a,b]\Longrightarrow \|x\|_{i'}\le c(\|a+b\|_{i'}+\|a-b\|_{i}).\end{equation}
The constant $c$ and the index $i$ do not depend on $a,b$.
\end{lem}
\begin{rem}In the sequel we do not use compactness of an interval. \end{rem}
{\it Proof of lemma \ref{dfg55}.}
Let us shift the set $[a,b]$ and consider a set
$$J=[a-s,b-s],\quad s=\frac{a+b}{2}.$$ The set $[a,b]$ is compact iff the set $J$ is compact.

Consider a projection $P_ny=\sum_{k=1}^ny_ke_k.$ Each set $$K_n=P_n(J)\subset J$$ is compact since it is a closed and bounded subset of $\mathbb{R}^n$.

Show that the sets $\{K_n\}$ form $\varepsilon$-nets in $J$.

Indeed, take any element $y\in J$ and present it as follows
$$y=P_ny+q_n,\quad q_n=\sum_{k=n+1}^\infty y_ke_k,\quad |y_k|\le r_k=\frac{b_k-a_k}{2}.$$

A series  $R_n=\sum_{k=n+1}^\infty r_ke_k\in E$
 is a tail of the expansion of the element $(b-a)/2$ and thus for all $i$ it follows that
$\|R_n\|_i\to 0.$
Observe that
$$q_n=\mathcal M_\lambda R_n,$$ where
$\lambda=\{\lambda_j\},\quad \lambda_j=y_j/r_j$ provided $r_j\ne 0$ and   $\lambda_j=0$  otherwise.

Theorem \ref{e5b5j} implies that for any $i'$ there exists $i$ such that
\begin{equation}\label{xdfgb0088}\|q_n\|_{i'}=\|\mathcal M_\lambda R_n\|_{i'}
\le c\|R_n\|_i\to 0.\end{equation}

The  limit  in the last part of formula (\ref{xdfgb0088}) is uniform in $y\in J$.
This proves the lemma in the part of compactness.

From the formulas above it also follows that $y=\mathcal M_\lambda R_0,\quad y\in J,$
$$\|y\|_{i'}\le c\|R_0\|_i=c\Big\|\frac{a-b}{2}\Big\|_i.$$
After the back shift we readily obtain estimate (\ref{dfgh77}).

Lemma \ref{dfg55} is proved.

\begin{lem}\label{dfh66pp}Let $ W\subset E$ be a nonempty set with an upper bound $\overline w\in E$:
$$w\in W\Longrightarrow w\ll \overline w.$$

Then the element
$\sup W$ is well defined and $\sup W\ll \overline w$.
\end{lem}
Indeed, the  assertion follows from lemma \ref{srg3ee}:  take any $x\in W$; then $\sup W\in[x,\overline w]$.

\begin{lem}\label{dsfg2ww}
Let a  function $F:I_T\times E\to \mathbb{R}$ be  left continuous in the second argument and a function $t\mapsto F(t,x)$ be integrable on $I_T$ for each $x\in E$.

Suppose also that $F$ is monotone: $$x\ll y\Longrightarrow F(t,x)\le F(t,y),\quad \forall t\in I_T.$$

Assume that a function $u:I_T\to E$,
$$u(t)=\sum_{k=1}^\infty u_k(t)e_k$$ is such that all the functions $u_k:I_T\to\mathbb{R}$ are Lebesgue measurable   and for some  $a,b\in E,\quad a\ll b$ and for   all $t\in I_T$ one has $$u(t)\in[a,b].$$

Then a mapping $t\mapsto F(t,u(t))$ is integrable on $I_T$.
\end{lem}
{\it Proof of lemma \ref{dsfg2ww}.}
From \cite{folland} we know that for each $k$ there exists a sequence  of simple functions $\ph_{k,j}(t)$ such that
$$a_k\le \ph_{k,j}\le \ph_{k,j+1}\le b_k,\quad a=\sum_{r=1}^\infty a_re_r,\quad b=\sum_{r=1}^\infty b_re_r$$ and $\ph_{k,j}\to u_k$ pointwise for each  $t\in I_T$ as $j\to\infty.$

Introduce the following functions
$$[a,b]\ni U_j(t)=\sum_{k=1}^j\ph_{k,j}(t)e_k+\sum_{r=j+1}^\infty a_re_r,\quad F_j(t)=F(t,U_j(t)).$$
Each function $U_j$ has a finite set of values in $E$ and $U_j\ll U_{j+1}.$ Thus $F_j:I_T\to \mathbb{R}$ is integrable.

Let us show that  $U_j\to u$ pointwise in $E$. Indeed, consider an estimate:
\begin{align}
\|U_j(t)&-u(t)\|_i\le \Big\|\sum_{k=1}^N(\ph_{k,j}(t)-u_k(t))e_k\Big\|_i\nonumber\\
&+\Big\|\sum_{k=N+1}^j(\ph_{k,j}(t)-u_k(t))e_k+\sum_{k=j+1}^\infty(a_k-u_k(t))e_k\Big\|_i,\quad j>N.\nonumber
\end{align}
The first summand in the right hand side of this formula vanishes as $j\to\infty$. By  lemma \ref{dfg55} the second summand is estimated from above in terms of
$$\Big\|\sum_{k=N+1}^ja_k e_k\Big\|_{i_1},\quad \Big\|\sum_{k=N+1}^ja_k e_k\Big\|_{i_2},\quad
\Big\|\sum_{k=N+1}^jb_k e_k\Big\|_{i_1},\quad \Big\|\sum_{k=N+1}^jb_k e_k\Big\|_{i_2}$$
and
$$\Big\|\sum_{k=j+1}^\infty a_k e_k\Big\|_{i_1},\quad \Big\|\sum_{k=j+1}^\infty a_k e_k\Big\|_{i_2},\quad
\Big\|\sum_{k=j+1}^\infty b_k e_k\Big\|_{i_1},\quad \Big\|\sum_{k=j+1}^\infty b_k e_k\Big\|_{i_2}.$$
These terms  vanish as $N\to \infty$.

 So that $F_j(t)\to F(t,u(t))$ pointwise. On the other hand
 $$F(t,a)\le F_j(t)\le F(t,b).$$
Therefore the assertion of the lemma follows from the Dominated convergence theorem.

Lemma \ref{dsfg2ww} is proved.

\begin{lem}\label{sfg00jj}
Take a function $u:I_T\to [a,b]\subset E,$
$$u(t)=\sum_{k=1}^\infty u_k(t)e_k$$ with  $u_k$  Lebesgue measurable. Then a function
$t\mapsto f(t,u(t))$ is integrable in $I_T$.\end{lem}{\it Proof of lemma \ref{sfg00jj}.}
Indeed, from lemma \ref{dsfg2ww} we know that the functions $f_k(\cdot,u(\cdot)):I_T\to\mathbb{R}$ are integrable and  $f_k(t,a)\le f_k(t, u(t))\le f_k(t,b)$.

Introduce  functions
$$\phi_n(\cdot)=\sum_{k=1}^nf_k(\cdot,u(\cdot))e_k$$ and observe that $\phi_n(\cdot)\to f(\cdot,u(\cdot))$ pointwise  in $E$.

The functions $\phi_n$ are integrable. Moreover,
$$ \mathcal M_{\lambda_n}f(t,a)\ll \phi_n(t)\ll \mathcal M_{\lambda_n}f(t,b),\quad
\lambda_n=(\underbrace{1,\ldots,1}_{\mbox{$n$ times}},0,0,\ldots) $$
and from
  lemma \ref{dfg55} and theorem \ref{e5b5j} it follows that for any $i\in\mathbb{N}$ there are $i',i''\in\mathbb{N}$ and a constant $c>0$ such that for  all $t$ one has
$$\|\phi_n(t)\|_{i}\le c(\|f(t,a)\|_{i'}+\|f(t,b)\|_{i'}+\|f(t,a)\|_{i''}+\|f(t,b)\|_{i''}).$$  By the statement of the problem the function in the right side of this inequality is integrable (for details see section \ref{xfhb77}).

The Dominated convergence theorem (theorem \ref{sdgaa} below) concludes the proof.

Lemma \ref{sfg00jj} is proved.

\subsection{Proof of Theorem \ref{dsfg5}}
\label{sxdfg55}
Introduce a set
\begin{align*}S=\Big\{u(t)=&\sum_{k=1}^\infty u_k(t)e_k\in E\,\Big|\,u_k\quad \mbox{are lower semicontinuous},\\&
x_*\ll u(t)\ll \Phi(u)(t),\quad t\in I_T\Big\},\end{align*}
where
$$\Phi(u)(t)=\hat x+\int_0^tf(\xi,u(\xi))d\xi.$$

The set $S$ is not empty: $x_*\in S$.

By formula (\ref{fgb00}) if $u\in S$ then one has
\begin{equation}\label{sfhhhhg00jj}\Phi(u)(t)\in[x_*,\hat x+T|C|],\quad u(t)\in[x_*,\hat x+T|C|],\quad t\in I_T.\end{equation}
Observe also that $\Phi(S)\subset S$.
Indeed, this follows by lemma \ref{sfg00jj} from the first inclusion of (\ref{sfhhhhg00jj}) and monotonicity of the mapping $f$.

The set $S$ is partially ordered by the following binary relation. For any $u,v\in S$ by definition put
$$u\prec v\Longleftrightarrow u(t)\ll v(t)\quad \forall t\in I_T.$$

Let $S(t)$ denote the following set $\{\eta(t)\mid \eta\in S\}\subset[x_*,\hat x+T|C|].$

By lemma \ref{dfh66pp} a function $\tilde x(t)=\sup S(t)$ is correctly defined. We then obtain
$$\tilde x(t)=\sum_{k=1}^\infty \tilde x_k(t)e_k\in[x_*,\hat x+T|C|].$$
From \cite{edvards} we know that all the functions $\tilde x_k$ are lower semicontinuous.

Furthermore for any $w\in S$ and for any $t\in I_T$ it follows that\begin{equation}\nonumber
\begin{split} w(t)&\ll \tilde x(t)\Longrightarrow f(t,w(t))\ll f(t,\tilde x(t))\nonumber\\
&\Longrightarrow \int_0^tf(\xi,w(\xi))d\xi\ll
\int_0^tf(\xi,\tilde x(\xi))d\xi.\nonumber
\end{split}
\end{equation}
And thus
$$w(t)\ll\hat x+\int_0^tf(\xi,w(\xi))d\xi\ll \hat x +\int_0^tf(\xi,\tilde x(\xi))d\xi.$$
The last estimate holds for all $w\in S$. This implies
$$\tilde x(t)\ll \hat x +\int_0^tf(\xi,\tilde x(\xi))d\xi$$ and $\tilde x\in S.$

By definition of $\tilde x$ one obtains $$\tilde x(\cdot)\prec\Phi(\tilde x(\cdot))\in S\Longrightarrow \tilde x(\cdot)=\Phi(\tilde x(\cdot))$$ and $\tilde x$ is the desired solution to problem (\ref{sdgf5500}).

The theorem is proved.

\subsection{Proof of Theorem \ref{sdfrg55oo}.}
Let $x\in Q$ be a solution. Then for all $t\in I_T$ we have
$$x(t)=\hat x+\int_0^t f(\xi,x(\xi))d\xi\ll \hat x+T|C|$$
and thus $\overline x(t)\ll \hat x+T|C|.$
Sinse $\tilde x\in Q$ we also obtain $x_*\ll\tilde x(t)\ll \overline x(t).$

Recall that all the functions $\overline x_k(t)$ from the expansion
$$\overline x(t)=\sum_{k=1}^\infty \overline x_k(t)e_k$$ are lower semicontinuous \cite{edvards}.

Therefore by lemma \ref{sfg00jj} the function $\xi\mapsto f(\xi,\overline x(\xi))$ is integrable and the following formula is correct:
$$x(t)=\hat x+\int_0^t f(\xi,x(\xi))d\xi\ll \hat x+\int_0^tf(\xi,\overline x(\xi))d\xi$$
and thus
$$ \overline x(t)\ll\hat x+\int_0^tf(\xi,\overline x(\xi))d\xi.$$

So that $\overline x\in S$ and then $\tilde x=\overline x$.

The theorem is proved.

\subsection{Proof of Theorem \ref{sxdfg1100}}
\label{xdfb7khkhkh}
Recall that a function $w:I_T\to\mathbb{R}$ is absolutely continuous iff it can be presented in the form
$$w(t)=\int_0^tp(\xi)d\xi,\quad p\in L^1(I_T).$$ In this case the classical derivative $\dot w$ exists for almost all $t$ and
$$\dot w(t)=p(t)\quad\mbox{almost everywhere}.$$
From lemma \ref{sfg00jj} we know that $f(\cdot,x(\cdot))\in L^1(I_T,E)$.
Thus
$$(\psi,x(t))=(\psi,\hat x)+\Big(\psi,\int_0^tf(\xi,x(\xi))d\xi\Big)=(\psi,\hat x)+\int_0^t\big(\psi,f(\xi,x(\xi))\big)d\xi$$ is an absolutely continuous function and
\begin{equation}\label{sdfg09gfv}
\frac{(\psi,x(t+h))-(\psi,x(t))}{h}=\Big(\psi,\frac{x(t+h)-x(t)}{h}\Big),\quad h\ne 0.\end{equation}
For almost all $t$ one has
\begin{equation}\label{srg5tt}\frac{d}{dt}(\psi,x(t))=\big(\psi,f(t,x(t))\big).\end{equation}

The space $E$ is separable so that $E''$ is also separable and thus $E'$ is separable.
Let $\Psi\subset E'$ be a countable dense set.

Since $x_*\ll x(t)$ we obtain
$f(t,x(t))\in [f(t,x_*),C]$ and  for any $i$ and for almost all $t$ we have
\begin{align}\Big\|&\frac{1}{h}\int_t^{t+h}f(\xi,x(\xi))d\xi\Big\|_i\le \frac{1}{|h|}\Big|\int_t^{t+h}\|f(\xi,x(\xi))\|_ids\Big|\nonumber\\
&\le \frac{1}{|h|}\Big|\int_t^{t+h}\tilde c_i\big(\|f(\xi,x_*)+C\|_i+\|f(\xi,x_*)-C\|_{i'}\big)d\xi\Big|
\le c_i(t).\label{dfg5000}\end{align}
Here $\tilde c_i$ is a positive constant.

Let $\Theta_\psi$ be a set of values $t\in I_T$ for which formula (\ref{srg5tt}) does not hold. The measure of $\Theta_\psi$ is equal to zero so that the measure of a set $$\Theta=\bigcup_{\psi\in \Psi} \Theta_\psi$$ equals zero as well.

By the same reason a set
$$\tilde\Theta=\Theta\cup\{\mbox{$t$ for which formula (\ref{dfg5000}) does not hold}\}$$ is of measure zero too.

Introduce a set $I_T'=I_T\backslash \tilde\Theta.$ Consider an element
$$X_{t,h}=\frac{x(t+h)-x(t)}{h}$$ as a linear function on $E'$.

For each $t\in I'$ one therefore obtains $\|X_{t,h}\|_i\le c_i(t).$

Fix $t\in I_T'$. From formula (\ref{sdfg09gfv}) for any $\psi\in\Psi$ one gets
$$(\psi,X_{t,h})\to \frac{d}{dt}(\psi,x(t)),\quad h\to 0.$$From the Banach-Steinhaus theorem \cite{edvards} it follows that
the limit $$(q(t),\psi)=\lim_{h\to 0}(\psi,X_{t,h})$$ exists for all $\psi\in E'$. For each $t\in I_T'$ the element  $q(t)\in E''$. By reflexivity of $E$ we can regard it as $q(t)\in E$ and put $\dot x=q$.

The theorem is proved.

\section{Appendix: A sketch  of the Lebesgue integration theory for functions with values in Fr\'echet spaces}\label{xfhb77}

Let $(M,\Sigma,\mu)$ be a measure space with $\sigma-$algebra $\Sigma$ and a measure $\mu:\Sigma\to [0,\infty]$;
and let $(F,\{\|\cdot\|_i\}_{i\in\mathbb{N}})$ be a Fr\'echet space.

\begin{df}We shall say that $f:M\to F$ is a step function if

1) $f$ has a finite set of values $\{a_1,\ldots,a_\nu\}\subset F$;

2) for each $k=1,\ldots,\nu$ a set $A_k=f^{-1}(a_k)$ is measurable;

3) $\mu\{x\in M\mid f(x)\ne 0\}<\infty.$

The integral of a step function is defined as follows
$$\int f=\sum_{a_k\ne 0} a_k\mu(A_k).$$
\end{df}
The integral possesses the standard elementary properties.

It is clear that the  step functions form a vector space which we denote by $\mathrm{St}(M,F)$. Being equipped with a collection of seminorms $$|f|_i=\int\|f\|_i, \quad \|f\|_i\in \mathrm{St}(M,\mathbb{R})$$ the space $\mathrm{St}(M,F)$ becomes a locally convex space.

\begin{df} We shall say that $\{f_n\}\subset \mathrm{St}(M,F)$ is an approximating sequence for a function $g:M\to F$ if

1) $\{f_n\}$ is an $L^1-$ Cauchy sequence:
for each $\eps>0$ and for each $s\in\mathbb{N}$ there exists a number $l$ such that
$$i,m>l\Longrightarrow |f_i-f_m|_s<\eps,$$
and

2) $f_n\to g$ almost everywhere.\end{df}
We use $\mathcal L$ to denote the space of functions $g$ that have an approximating sequence. The functions from $\mathcal L$ we call integrable functions.

\begin{theo}\label{dfg600oo}
Let $\{f_n\}\subset \mathrm{St}(M,F)$ be an $L^1-$Cauchy sequence. Then it contains a subsequence $\{f_{j_p}\}$  that converges almost everywhere to a function $g:M\to F$ and for any $\eps>0$ there exists a measurable set $P,\quad \mu(P)<\eps$
such that $f_{j_p}$ converges to $g$ uniformly in $M\backslash P$ i. e.
for any $s\in\mathbb{N}$ one has
$$\lim_{j_p\to\infty} \sup_{x\in M\backslash P}\|f_{j_p}(x)-g(x)\|_s=0.$$\end{theo}
{\it Proof of theorem \ref{dfg600oo}.}
The idea of lemma 3.1 \cite{lang} is not destroyed if we  substitute  a seminorm $\|\cdot\|_s$ instead of the norm in its proof.  Consequently from the proof of this lemma we conclude that the sequence
$\{f_n\}$ contains a  subsequence $f_n^1$ such that for almost all $x\in M$ one has
$$\|f_n^1(x)-f_m^1(x)\|_1\to 0\quad \mbox{as}\quad n,m\to\infty$$ and for some measurable $P_1,\quad \mu(P_1)<\eps/2$ it follows that
$$\lim_{j,i\to\infty} \sup_{x\in M\backslash P_1}\|f^1_j(x)-f^1_i(x)\|_1=0.$$

Repeating the argument we extract a subsequence $\{f_n^2\}\subset\{f_n^1\}$ such that
for almost all $x\in M$ one has
$$\|f_n^2(x)-f_m^2(x)\|_2\to 0\quad \mbox{as}\quad n,m\to\infty$$ and for some measurable $P_2,\quad \mu(P_2)<\eps/2^2$ it follows that
$$\lim_{j,i\to\infty} \sup_{x\in M\backslash P_2}\|f^2_j(x)-f^2_i(x)\|_2=0$$
and so forth.

Thus for almost all $x\in M$ the diagonal sequence $f_{j_p}(x):=f^p_p(x)$ is a Cauchy sequence in $F$ and therefore it is convergent to some function $g:M\to F$ almost everywhere. It is also clear that
$$P=\bigcup_{i=1}^\infty P_i,\quad \mu(P)<\eps.$$

The theorem is proved.

\begin{theo}\label{dfh66cfr}Let $\{f_n\},\{h_n\}\subset \mathrm{St}(M,F)$ be $L^1-$Cauchy sequences that converge almost everywhere to the same function $f$.

Then one has $$\lim_{n\to\infty}\int f_n=\lim_{n\to\infty}\int h_n.$$\end{theo}
The existence of the
limits above is obvious: for any $s\in\mathbb{N}$ we have
$$\Big\|\int f_n-\int f_m\Big\|_s\le \int\|f_n-f_m\|_s\to 0\quad \mbox{as}\quad m,n\to\infty.$$
Thus $\int f_n$ is a Cauchy sequence in $F$.

{\it Proof of theorem \ref{dfh66cfr}.} The argument of lemma 3.2 \cite{lang} is adapted in the same way:
we carry out all the steps of the proof replacing the norm with an arbitrary seminorm from the  family $\{\|\cdot\|_s\}$.
This gives
$$|f_n-h_n|_s\to 0\quad \mbox{as $n\to\infty$},\quad \forall s\in \mathbb{N}.$$
The theorem is proved.

The last theorem allows to define the integral of a function $f\in\mathcal L.$ Indeed, let $f_n$ be an approximating sequence for $f$. Then by definition put
$$\int f=\lim_{n\to\infty}\int f_n.$$

Observe also  that for each $s$ we have $$\big|\|f_n\|_s- \|f\|_s\big|\le \|f_n-f\|_s\to 0$$ almost everywhere
and $$\int\big|\|f_n\|_s- \|f_m\|_s\big|\le \int\|f_n-f_m\|_s\to 0$$ as $m,n\to\infty$. So that  $\|f\|_s$ is an integrable function in the standard Lebesgue sense ($\|f_n\|_s$ is an approximating sequence for $\|f\|_s$).

\begin{lem}\label{sdfg6900}If $\{f_n\}$ is an approximating sequence for $f$ then it converges to $f$ relative the seminorms $|\cdot|_s$ i. e. for any $s\in\mathbb{N}$ it follows that
$|f-f_n|_s\to 0.$\end{lem}
Indeed, fix $n$; a sequence $h^s_{in}=\|f_i-f_n\|_s$ is an approximating sequence for $\|f-f_n\|_s$:
$$\int |h^s_{in}-h^s_{jn}|\le \int\|f_i-f_j\|_s;$$
$$h^s_{in}\to \|f-f_n\|_s\quad \mbox{almost everywhere as $i\to\infty$}$$
and $$\int \|f-f_n\|_s=\lim_{i\to\infty}\int h^s_{in}.$$
By the definition of an approximating sequence, the integral  $\int h^s_{in}$ vanishes as $i,n\to\infty$.

\begin{theo}\label{dfg05566hh}
The space $\mathcal L$ is complete relative the seminormes $|\cdot|_s$.\end{theo}
{\it Proof of theorem \ref{dfg05566hh}.}
Let $\{g_n\}\subset \mathcal L$ be an $L^1-$Cauchy sequence. There exists a sequence $\{f_n\}\subset \mathrm{St}(M,F)$
such that
$$\max_{1\le i\le n}|f_n-g_n|_i<\frac{1}{n}.$$
The sequence $f_n$ is an $L^1-$Cauchy sequence: for any $s$ it follows that
$$|f_j-f_l|_s\le |f_j-g_j|_s+|f_l-g_l|_s+|g_j-g_l|_s.$$
From theorem \ref{dfg600oo} we see that there is a subsequence $\{f_{n_l}\}\subset \{f_n\}$ such that $f_{n_l}\to g$ almost everywhere.

The sequence $\{f_{n_l}\}$ is an approximating sequence for $g$. By lemma \ref{sdfg6900} it follows that $|f_{n_l}-g|_s\to 0$.

To conclude the proof observe that
$$|g_n-g|_s\le |g_n-f_n|_s+|f_{n_l}-f_n|_s+|f_{n_l}-g|_s.$$

\begin{theo}\label{dxfgww}
Suppose that  a sequence $\{f_n\}\subset \mathcal L$ converges to an element $f\in\mathcal L$ relative the seminorms $\{|\cdot|_s\}$ . Then it contains a subsequence $\{f_{j_p}\}$  that converges almost everywhere to  $f$ and for any $\eps>0$ there exists a measurable set $P,\quad \mu(P)<\eps$
such that $f_{j_p}$ converges to $f$ uniformly in $M\backslash P$.\end{theo}
The proof of this theorem repeats the proof of theorem \ref{dfg600oo}.

Theorem \ref{dxfgww} implies in particular that for $f\in\mathcal L$ the following two facts are equivalent:
$$(f=0\quad \mbox{almost everywhere})\Longleftrightarrow (|f|_s=0\quad\forall s\in\mathbb{N}).$$

Combining theorems \ref{dxfgww} and \ref {dfg05566hh} we obtain
\begin{theo}\label{xfb0011}

Let a sequence $\{f_n\}\subset \mathcal L$ be an $L^1-$Cauchy sequence that converges almost everywhere to a function $f$. Then $f\in\mathcal L$ and $|f_n-f|_s\to 0$ for any $s$.\end{theo}

\begin{theo}[Dominated convergence theorem]\label{sdgaa}
Assume that a sequence $\{f_n\}\subset\mathcal L$ converges almost everywhere to a function $f$. Assume also that there exists a sequence of Lebesgue integrable functions $g_s:M\to\mathbb{R}$ such that the inequalities
$$\|f_n\|_s\le g_s$$ hold for all $n,s\in\mathbb N.$
Then  $f\in\mathcal L$ and $\{f_n\}$ converges to $f$ relative the seminorms $|\cdot|_s$.\end{theo}
{\it Proof of theorem \ref{sdgaa}.}
Indeed, introduce functions $$\psi^s_n(x)=\sup_{i,j\ge n}\|f_i(x)-f_j(x)\|_s.$$
For almost all $x$ we have
$$\psi^s_n(x)\to 0$$  as $n\to\infty$ and  $0\le \psi^s_n(x)\le 2g_s(x)$. By the standard finite dimensional Lebesgue integration theory these functions are measurable.

From the finite dimensional Dominated convergence theorem we obtain $\int\psi_n^s\to 0$ for each $s\in\mathbb{N}.$
Therefore $\{f_n\}$ is an $L^1-$Cauchy sequence and theorem \ref{xfb0011} concludes the proof.

\end{document}